\documentclass[aap,MSNbibl,dvips]{arximspdf}
\usepackage{graphicx}

\doi{10.1214/12-AAP877} 
\volume{23}
\issue{4}
\pubyear{2013}
\firstpage{1469}
\lastpage{1493}

\makeatletter
\newcommand{\rrvert}{\vert}
\newcommand{\llvert}{\vert}

\newcommand{\R}{\mathbb{R}}

\newcommand{\N}{\mathbb{N}}
\newcommand{\E}{\mathbb{E}}
\renewcommand{\P}{\mathbb{P}}
\renewcommand{\R}{\mathbb{R}}
\renewcommand{\d}{{\mathrm d}}
\renewcommand{\t}{{\mathcal T}}

\newtheorem{theorem}{Theorem}
\newtheorem{proposition}{Proposition}
\newtheorem{lemma}{Lemma}
\newtheorem{corollary}{Corollary}

\makeatother

\begin{document}
\begin{frontmatter}

\title{The cut-tree of large Galton--Watson trees
and the~Brownian CRT}
\runtitle{Cut-tree}

\begin{aug}
\author[A]{\fnms{Jean} \snm{Bertoin}\corref{}\ead[label=e1]{jean.bertoin@math.uzh.ch}}
\and
\author[B]{\fnms{Gr\'egory} \snm{Miermont}\thanksref{t1}\ead[label=e2]{Gregory.Miermont@math.u-psud.fr}}
\runauthor{J. Bertoin and G. Miermont}
\affiliation{Universit\"at Z\"urich and Universit\'e Paris-Sud}
\address[A]{Institut f\"ur Mathematik\\
Universit\"at Z\"urich\\
Winterthurerstrasse 190\\
CH-8057 Z\"urich\\
Switzerland\\
\printead{e1}} 
\address[B]{Equipe de Probabilit\'es, Statistiques\\
\quad et Mod\'elisation\\
Universit\'e Paris-Sud\\
B\^{a}timent 425\\
91405 Orsay Cedex\\
France\\
\printead{e2}}
\end{aug}

\thankstext{t1}{Supported by ANR Grant ANR-08-BLAN-0190 (A3).}

\received{\smonth{1} \syear{2012}}
\revised{\smonth{5} \syear{2012}}

%
\begin{abstract}
Consider the edge-deletion process in which the edges of some finite
tree $T$ are removed one after the other in the uniform random order.
Roughly speaking, the cut-tree then describes the genealogy of
connected components appearing in this edge-deletion process. Our main
result shows that after a proper rescaling, the cut-tree of a critical
Galton--Watson tree with finite variance and conditioned to have size
$n$, converges as $n\to\infty$ to a Brownian continuum random tree
(CRT) in the weak sense induced by the Gromov--Prokhorov topology. This
yields a multi-dimensional extension of a limit theorem due to Janson
[\textit{Random Structures Algorithms} \textbf{29} (2006) 139--179] for
the number of random cuts needed to isolate the root in Galton--Watson
trees conditioned by their sizes, and also generalizes a recent result
[\textit{Ann. Inst. Henri Poincar\'e Probab. Stat.} (2012)
\textbf{48} 909--921]
obtained in the special case of Cayley trees.
\end{abstract}

%

\end{frontmatter}

\section{Introduction and main results}\label{sec1}
\subsection{Motivations}\label{sec11}

Random destruction of combinatorial trees is an old topic which can be
traced back more than 40 years ago to
the work of Meir and Moon~\cite{MM}.
Let $T$ be a rooted tree on a finite set of vertices. Imagine that we pick
a vertex uniformly at random and destroy it together with the entire
subtree generated by that vertex. We iterate in an obvious way until
the root is picked and are interested in the
number $N(T)$ of steps of this algorithm.

The present paper has been motivated by the following result due
to Janson~\cite{Janson} who treated the case where $T$ is a large
Galton--Watson tree. More precisely, consider the genealogical tree of
a branching process having a critical reproduction law with finite
variance $\sigma^2>0$, and let $\t_n$ be a version of this tree
conditioned to have exactly $n$ vertices, assuming implicitly that the
probability of that event is positive. Then Janson established that
$N(\t_n)/(\sigma\sqrt n)$ converges weakly as $n\to\infty$ to the
Rayleigh distribution which has density $x\exp(-x^2/2)$ on $\R_+$.
See also Panholzer~\cite{Panholzer} for the same result in a less
general setting, and Abraham and Delmas~\cite{AD} for a recent
contribution and further references.

The following extension has been recently obtained in~\cite{Be}. Let
$T_n$ be a uniform Cayley tree with $n$ vertices; it is well known that
this corresponds to a special case of conditioned Galton--Watson trees,
namely, when the reproduction law is Poisson. Given $T_n$, distinguish
$k$ vertices uniformly at random, where $k$ is some fixed integer. Then
remove an edge uniformly at random and independently of the
distinguished vertices. This disconnects $T_n$ into two subtrees. If
one of these subtrees does not contain any of the distinguished
vertices, then we destroy it entirely, else we keep the two subtrees.
We iterate until each and every distinguished vertex has been isolated
and denote by $Y(T_n,k)$ the number of steps.
Then, according to Lemma 1 in~\cite{Be}, $Y(T_n,k)/ \sqrt n$ converges
weakly as $n\to\infty$ to the Chi distribution with parameter $2k$,
which has density
\[
\frac{2^{1-k}}{(k-1)!}x^{2k-1} \exp\bigl(-x^2/2\bigr)
\]
on $\R_+$.
This result has also been very recently recovered by~\cite{ABH} using
a different approach.

The Chi$(2k)$ distribution occurs as the law of the length $L_k(\mathbf
{T})$ of a Brownian continuum random tree (CRT) $\mathbf{T}$ reduced
to $k$ leaves picked uniformly at random, as can be seen from Aldous
\cite{ACRTIII}, Lemma 21. The appearance of the Brownian CRT in this
framework should not come as a surprise since it is well known that if
we assign length $1/\sqrt n$ to each edge of $T_n$, then the latter
converges weakly to a Brownian CRT $\mathbf{T}$ as $n\to\infty$. We
stress, however, that the rescaled Cayley tree $n^{-1/2}T_n$ and
$n^{-1/2} Y(T_n,k)$ \textit{do not} converge jointly in distribution
toward $\mathbf{T}$ and $L_k(\mathbf{T})$.

The proof in~\cite{Be} of the extension above of Janson's result
relies on three crucial features. First, the observation due to Pitman
\cite{Pi} that random deletion of edges in a uniform Cayley tree
yields a remarkable fragmentation process; second, a general limit
theorem due to Haas and Miermont~\cite{HM} for so-called branching
Markov trees; third, the characterization of the Brownian
fragmentation in~\cite{Be1}. More precisely, the fragmentation process
that results from the repeated deletion of edges in a uniform Cayley
tree can be represented by a Markov branching tree whose law is
explicitly known. In this setting, $Y(\t_n,k)$ corresponds to the
length of
this Markov branching tree reduced to $k$ leaves picked uniformly at
random. Thanks to the limit theorem of Haas and Miermont, one then
checks that this Markov branching tree with lengths rescaled by a
factor $1/\sqrt n$ converges weakly, and the limit can then be
identified as another Brownian CRT, say ${\mathbf T}'$, using the
characterization of the fragmentation process at heights induced by
the latter. As a consequence, $n^{-1/2} Y(\t_n,k)$ converges weakly
to $L_k(\mathbf{T'})$ and hence, to the Chi$(2k)$ law. Unfortunately,
this approach only works for Cayley trees; as for other conditioned
Galton--Watson trees, random edge deletion does not yield, in general, a
Markov branching tree as above, and the entire structure of the proof
collapses. Nonetheless, the fact that Janson's result is valid for any
critical Galton--Watson tree with finite variance suggests that the
same should also hold for its natural extension to $k$ vertices for
$k\geq2$.

The first purpose of this work is to show that this is indeed the case. For
the sake of convenience, we shall deal with a slightly modified model
in which we distinguish edges rather than vertices, and
which is easily seen to have the same asymptotic behavior as the
former. The precise framework and result are presented in Section \ref
{sec12} below.

Our main goal, in the spirit of~\cite{Be}, will be to prove a
convergence result
for the genealogy induced by the edge-deletion
procedure, even though this process does not satisfy, in general, the Markov
branching property of~\cite{HM}. In Section~\ref{sec13}, we introduce the
cut-tree of a finite tree, which roughly speaking records the genealogy
of blocks in the edge-deletion process which consists of removing edges
of that tree one after the other and in uniform random order. In
Section~\ref{sec14}, we define the cut-tree of a Brownian CRT $\mathbf{T}$,
relying on a Poissonian logging process on the skeleton of $\mathbf
{T}$ which has been constructed by Aldous and Pitman~\cite{AP} to
study the so-called standard additive coalescent. Our main result
claims the joint weak convergence of $\t_n$ and its cut-tree suitably
rescaled toward their continuous counterparts, namely, $\mathbf{T}$ and
$\operatorname{cut}(\mathbf{T})$; it is stated in Section~\ref{sec15}.

After this long Introduction, the rest of this paper will be organized
as follows. Section~\ref{sec2} is devoted to preliminary results that
will be
used in the proof of Theorem~\ref{T}, and the latter is established in
Section~\ref{sec3}. Section~\ref{sec4} is devoted to the proof of a
technical bound,
relying partly on an invariance property under random re-planting for
Galton--Watson trees, which may be of independent interest.

\subsection{The number of cuts needed to find a few edges}\label{sec12}

\begin{figure}

\includegraphics{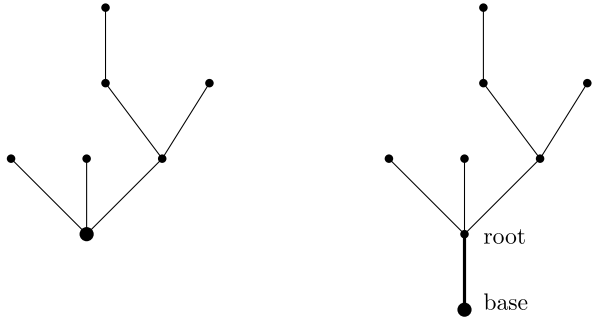}

\caption{Planting.}
\label{figplanting}
\end{figure}

It will be convenient in the sequel to work with a slight modification
of the trees under consideration. Consider a (rooted) tree $T$ on a
set of $n$ vertices, say $[n]=\{1,\ldots, n\}$; we add a new vertex
which we call the \textit{base} and link it to the root of $T$ by a new
edge. This gives a \textit{planted} tree which we denote by $\bar{T}$.
See Figure~\ref{figplanting}.

%
%

The set $\bar{E}$ of edges of $\bar{T}$ is thus given by the set
${E}$ of edges of $T$ plus the new edge connecting the base to the root.
We consider $\bar{E}$ as a set of vertices, and endow it with a
natural tree structure
by declaring that $e$ and $e'$ are neighbors in $\bar{E}$
if and only if they are adjacent in $\bar{T}$. Plainly, this yields a
tree which is isomorphic to $T$;
more precisely, the map $v\dvtx\bar{E} \to[n]$ that associates to an
edge $e$ of $\bar{T}$,
its extremity $v(e)\in[n]$ which is the farthest
away from the base vertex in $\bar{T}$,
is bijective and
preserves the tree structures.
Any statement expressed in terms of the edges of $\bar{T}$ can thus be
rephrased in terms of the vertices of $T$ and vice versa. For a
technical reason, it will be slightly simpler for us to work with the
edge-version rather than the vertex-version of conditioned
Galton--Watson trees.

As before, we consider a critical reproduction law $\nu$ with finite
variance $\sigma^2>0$. Denote by $p$ the greatest common divisor of
the support of $\nu$; and observe that when $n$ is a sufficiently
large integer, the total population of a
Galton--Watson process with reproduction law $\nu$ generated by a
single ancestor equals $n$ with positive probability if and only if
$n-1\in p\N$.

Let $\t_n$ denote a version of a Galton--Watson tree with reproduction
law $\nu$ conditioned to have exactly $n$ vertices which are
enumerated, for instance, in the breadth-first search order to yield a
tree-structure on $[n]$ as required. Of course, the vertex
corresponding to the ancestor serves as the root. We shall implicitly
restrict our attention to the case $n-1\in p\N$, so this conditioning
makes sense provided that $n$ is large enough. Recall that the
associated planted tree is denoted by $\bar{\t}_n$ and has $n$ edges.

Next, for every integer $k\geq1$, given $\bar{\t}_n$, we distinguish
$k$ edges in $\bar{\t}_n$ uniformly at random. Conditionally on
$\bar{\t}_n$, we pick an edge uniformly at random and independently of
these $k$ distinguished edges. We remove\vspace*{1pt} it from $\bar{\t}_n$; this
disconnects $\bar{\t}_n$ into two subtrees. We then only consider
subtrees that contain at least one of the distinguished edges,
discarding, if necessary, those that contain no distinguished edge. We
iterate until each and every distinguished edge has been removed, and
write $N(\bar{\t}_n,k)$ for the number of steps. We shall prove the
following.

%
%
\begin{proposition}\label{P1}
In the notation above,
\[
\frac{1}{\sigma\sqrt n} N(\bar{\t}_n,k)
\]
converges in distribution as $n\to\infty$ to the Chi distribution
with parameter $2k$, that is,
\[
\frac{2^{1-k}}{(k-1)!}x^{2k-1} \exp\bigl(-x^2/2\bigr) \,\d x,\qquad x>0.
\]
\end{proposition}

%
%
\begin{figure}[b]

\includegraphics{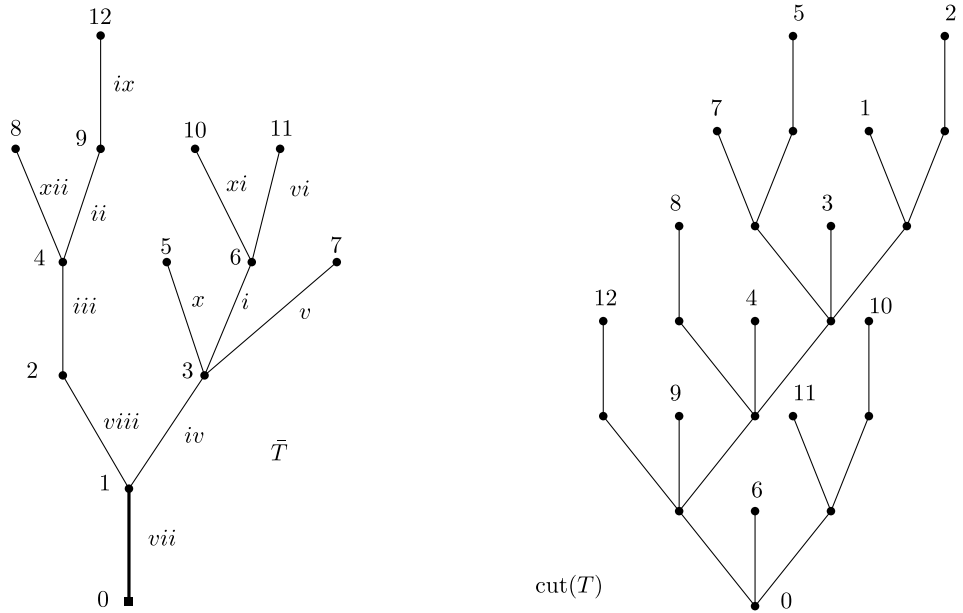}

\caption{The tree $\operatorname{cut}(T)$ of a planted tree $\bar T$. The
vertices are labeled in Arabic numerals in breadth-first order for this planar
representation of $\bar{T}$, while the order of deletion of the
edges is indicated in lowercase Roman numerals. Each internal vertex in
the tree
to the right is naturally labeled by the set of leaves that lie in
the subtree above this node.}
\label{figcut}
\end{figure}

\subsection{Cut-trees of finite trees}\label{sec13}

We can be more accurate by keeping track of the genealogy induced by
the edge-deletion process depicted above. More specifically, let $T$ be
a rooted tree with $n$ vertices and ${\bar T}$
its planted version. Recall that ${\bar T}$ has $n$ edges which are
naturally enumerated by the map $v\dvtx\bar{E}\to[n]$ described in the
preceding section.

Then we entirely destroy $\bar T$ by inductively removing its edges,
uniformly at random one after the other.
For $j=1, \ldots, n$, we denote by $i_j$ the label of the edge that is
removed at the $j$th step,
so $(i_1, \ldots, i_n)$ is a uniform random permutation of~$[n]$.
We partly encode this edge-deletion process by another tree, which we
denote by
$\operatorname{cut}(T)$ and construct as follows
(see Figure~\ref{figcut}).
For every $r=0, \ldots, n-1$, let
$\Pi(r)$ be the partition of $E(r):=\{1,2,\ldots,n\}\setminus\{i_1,
\ldots, i_r\}$
obtained by specifying that two elements $j$ and $j'$ in $E(r)$ are in
the same block
of $\Pi(r)$ if and only if either $j=j'$ or the edges with labels $j$
and $j'$ are still connected in the forest obtained from $\bar T$ by
deletion of the $r$ first edges with labels $i_1, \ldots, i_r$.
The family of the blocks (without repetition) of the partitions $\Pi
(r)$ for $r=0,\ldots, n-1$
forms the set of internal nodes of $\operatorname{cut}(T)$, the initial
block $[n]$ of $\Pi(0)$ being seen as the root. The leaves of $
\operatorname{cut}(T)$ are given by $1, \ldots, n$; we stress that a
singleton $\{
i\}$ may appear as an internal node of $\operatorname{cut}(T)$ and should
not be confused with the leaf $i$.

Now consider the $r$th step at which the edge labeled $i_r$ is removed,
and let
$B$ denote the block of $\Pi(r-1)$ which contains $i_r$. There are
three possibilities.
First, $B$ is reduced to the singleton $\{i_r\}$; in that case we draw
a single edge between the internal node $B=\{i_r\}$ and the leaf $i_r$.
Second, $B$ is not a singleton and $B\setminus\{i_r\}=B'$ is a block
of $\Pi(r)$; then we draw an edge between the internal nodes $B$ and
$B'$, and another edge between $B$ and the leaf $i_r$. Third, there are
two distinct blocks $B'$ and $B''$ of $\Pi(r)$ which result from $B$,
that is, $B=B'\sqcup B''\sqcup\{i_r\}$. Then we draw an edge between
the internal nodes $B$ and $B'$,
a second edge between the internal nodes $B$ and $B''$ and a third edge
between $B$ and the leaf $i_r$. If $T$ is a random
tree, we define $\operatorname{cut}(T)$ by first conditioning on $T$ and
then performing the above construction.

The main purpose of this work is to determine the asymptotic behavior
(in distribution)
of $\operatorname{cut}(\t_n)$ for a Galton--Watson tree $\t_n$ with size
$n$, as $n\to\infty$.
In this direction, it is appropriate to work in the framework of
pointed metric measure spaces.
More precisely, a finite tree $T$ with $n$ vertices can be identified as
$([n], d_n, \mu_n)$ where $[n]=\{1,\ldots, n\}$ is the set of
vertices, $d_n$ is the graph-distance on
$[n]$ induced by $T$ and $\mu_n$ the uniform probability measure on $[n]$.
We further retain the fact that $T$ is rooted by distinguishing in
$[n]$ the root vertex (usually~$1$, e.g., when $T$ is a
genealogical tree and vertices are labeled according to the
breadth-first search order).

It will be convenient to adopt a slightly different point of view for
$\operatorname{cut}(T)$,
by focusing on leaves rather than internal nodes.
More precisely, we set $[n]^0=[n]\cup\{0\}=\{0,1, \ldots, n\}$ where
$0$ corresponds to the root $[n]$ of $\operatorname{cut}(T)$ and
$1,\ldots, n$ to the leaves,
and consider the (random) metric measure space
$([n]^0, \delta_n, \mu_n)$ where $\delta_n$ is the (random) graph
distance on $[n]^0$ induced by $\operatorname{cut}(T)$, and, by a slight
abuse of notation, $\mu_n$ the uniform probability measure on $[n]$
extended by $\mu_n(0)=0$. It is easily seen that $\operatorname{cut}(T)$,
and in particular its combinatorial structure, can be recovered from
$([n]^0, \delta_n, \mu_n)$.

We stress that in this framework, the root of the cut-tree (which
corresponds to the additional point $0$ in $[n]^0$) has a crucial role.
Indeed, the height (distance to the root) of
the leaf $i$ in $\operatorname{cut}(T)$ is precisely the total number of
edge-removals from the successive blocks containing $i$ until the edge
$i$ is finally removed. More generally, the number of internal nodes of
the tree $\operatorname{cut}(T)$ reduced to its root and $k$ leaves, say
$\ell_1, \ldots, \ell_k$, coincides with the total number of
edge-removals from the blocks which contain at least one of those $k$
leaves until each and every one of the edges with labels $\ell_j$ in
$\bar T$ have been removed. Note also that this number differs from the
total length of that reduced tree by, at most, $k$ units. For this
reason, it will be important to recall that $0$ has been singled out in
the metric space $[n]^0$.

\subsection{Cut-tree of a Brownian CRT}\label{sec14}\label
{seccut-tree-continuum}

Aldous and Pitman~\cite{AP} considered a cutting process on the
Brownian CRT which bears obvious similarities with that defined in the
preceding section for finite trees. We recall and develop some
features in this setting that will be useful for our purpose.

Let us first recall some basic facts about topologies on metric measure
spaces that we will need. A pointed metric measure space is a
quadruple $(X,d,\mu,x)$ where $(X,d)$ is a complete metric space,
$x\in X$ and $\mu$ is a Borel probability measure on $(X,d)$.

Two such spaces $(X,d,\mu,x)$ and $(X',d',\mu',x')$ are called
isometry-equivalent if there exists an isometry
$f\dvtx\operatorname{supp}(\mu)\cup\{x\}\to X'$ (here $\operatorname
{supp}$ is the
topological support) such that $f(x)=x'$ and the image of $\mu$ by $f$
is $\mu'$. This defines an equivalence relation between pointed
metric measure spaces, and we note that the representatives
$(X,d,\mu,x)$ of a given isometry-equivalence class can always be
assumed to have $\operatorname{supp}(\mu)\cup\{x\}=X$.

The set $\mathbb{M}$ of (isometry-equivalence classes of) pointed
metric measure spaces is a Polish space when endowed with the
so-called Gromov--Prokhorov topology. Gromov's book~\cite{gromov99} and
the paper~\cite{GrPfWi} are good
references, although they deal with nonpointed spaces, which differs
from our setting only in a minor way. Recall that a sequence
$(X_n,d_n,\mu_n, x_n)$ of pointed measure metric spaces converges in
the Gromov--Prokhorov sense to $(X_\infty,d_\infty,\mu_\infty,
x_{\infty})$ if and only if the following holds: for $n\in
\N\cup\{\infty\}$, set $\xi_n(0)=x_n$ and let
$\xi_n(1),\xi_n(2),\ldots$ be a sequence of i.i.d. random variables
with law $\mu_n$, then the vector $(d_n(\xi_n(i),\xi_n(j)),0\leq
i,j\leq k)$ converges in distribution to
$(d_\infty(\xi_\infty(i),\xi_\infty(j)),0\leq i,j\leq k)$ for every
$k\geq1$.

Now recall that $\mathbf{T}$ denotes a Brownian CRT. It is endowed
with the uniform probability ``mass'' measure $\mu$ and the usual
distance $d$, and also comes with a distinguished point
called the root~\cite{ACRTIII}. Therefore,
$\mathbf{T}$ is viewed as a random variable in $\mathbb{M}$. Note that
the root plays the same role as a $\mu$-randomly chosen point in
$\mathbf{T}$, which is usually called the invariance property of
$\mathbf{T}$ under re-rooting.

The distance $d$ induces an extra length
measure $\lambda$, which is the unique $\sigma$-finite measure
assigning measure $d(x,y)$ to the geodesic path between $x$ and $y$ in
$\mathbf{T}$. Roughly speaking, the probability measure $\mu$ is
carried by the subset of leaves of $\mathbf{T}$ while the length
measure $\lambda$ rather lives on the skeleton, that is, the complement
of the set of leaves.

Conditionally on $\mathbf{T}$, we introduce the family $(t_i,
x_i)_{i\in I}$ of the atoms of a Poisson random
measure with intensity $\d t\otimes\d\lambda$, where $I$ is a
countable index set. We view these atoms as marks that are deposited
along the skeleton of $\mathbf{T}$ as time grows. Let $\mathbf{T}(t)$
be the ``forest'' obtained by removing the points $\{x_i\dvtx i\in
I,t_i\leq t\}$ that are marked before time $t$. For every $x\in\mathbf
{T}$ we
let $\mu_x(t)$ be the $\mu$-mass of the component of $\mathbf{T}(t)$ that
contains $x$, where by convention we let $\mu_x(t)=0$ if $x=x_i$ for
some $i\in I$ with $t_i\leq t$.

Aldous and Pitman~\cite{AP}
have observed that if $\xi$ denotes a random point in $\mathbf{T}$
distributed according to $\mu$ and independent of the Poisson point
process of marks on the skeleton, then the processes
%
%
\begin{equation}
\label{eq1} \bigl(\mu_{\xi}(t), t\geq0\bigr) \quad\mbox{and}\quad \bigl(1/
\bigl(1+\sigma(t)\bigr), t\geq0\bigr) \qquad\mbox{have the same law,}
\end{equation}
where $(\sigma(t), t\geq0)$ denotes the first-passage time process of
a linear Brownian motion.
Specializing results of~\cite{Be2} in this setting, we easily deduce
that if we define
\[
h_x=\int_0^\infty
\mu_x(t)\,\d t,\qquad x\in\mathbf{T},
\]
then
%
%
\begin{equation}
\label{eq2} h_{\xi} \mbox{ has the Rayleigh distribution}
\end{equation}
(see Section~\ref{sec32} below for details). As a consequence,
$0<h_x<\infty$
a.s. for $\mu$-almost every $x$.


We next add to $\mathbf{T}$ an extra point, denoted for simplicity by $0$,
and write $\mathbf{T}^0=\mathbf{T}\cup\{0\}$; $0$ serves, of course,
as distinguished element
of $\mathbf{T}^0$. We define a (random) function
$\delta$ of two arguments in $\mathbf{T}^0$ by setting
\[
\delta(0,0)=0,\qquad \delta(0, x)= \delta(x,0) = h_x\quad
\mbox{and}\quad
\delta(x,y) = \int_{t_{xy}}^\infty\bigl(
\mu_x(t)+\mu_y(t)\bigr)\,\d t,
\]
where for $x, y\in\mathbf{T}$ with $x\neq y$, $t_{xy}$ denotes
the (a.s. finite) smallest time $t$ when $x$ and $y$ belong to two
distinct components of $\mathbf{T}(t)$. Note that $t_{xy}$ is the
first time where a mark appears on the geodesic from $x$ to $y$, and
as such it has an exponential distribution of parameter $d(x,y)$.
Observe also that $\mu(\{x\in\mathbf{T}\dvtx\delta(0,x)=0\})=0$ a.s.
since $h_{\xi}>0$ a.s.
and that for every $y\in\mathbf{T}$, $\mu(\{x\in\mathbf{T}\dvtx
\delta(x,y)=0\})=0$ a.s. since $\mu_{\xi}(t)>0$ for all $t\geq0$,
a.s.

Let $\xi(0)=0$
and $(\xi(i),i\geq1)$ be an i.i.d. sequence with law $\mu$
conditionally given $\mathbf{T}$. We will see in Lemma~\ref{lem5}
below that the two random semi-infinite matrices
%
%
\begin{equation}
\label{eqa1} \bigl(d\bigl(\xi(i+1),\xi(j+1)\bigr)\dvtx i,j\geq0\bigr)
\quad\mbox
{and}\quad
\bigl(\delta\bigl(\xi(i),\xi(j)\bigr)\dvtx i,j\geq0\bigr)
\end{equation}
have the same distribution.\footnote{The shift of indices in the left-hand
side comes from the fact that the distinguished point $0$ is not the
root of $\mathbf{T}$, and formally it is not even an element of the
latter.} In particular, $\delta$ is a.s. a distance on the set
$\{\xi(i),i\geq0\}$, and $\mathcal{R}(k)=(\{\xi(i),0\leq i\leq
k\},\delta)$ can be understood as a \textit{consistent family} of random
rooted trees with, respectively, $k$ leaves in the sense of Aldous
\cite{ACRTIII}; here, the spaces $\mathcal{R}(k)$ have $k+1$
elements while they should really be ``trees with edge-lengths'' in
the context of Aldous' paper, but this is only a minor difference that
does not affect our discussion. Since $\mathbf{T}$ satisfies the
so-called \textit{leaf-tight property}
\[
\inf_{i\geq2}d\bigl(\xi(1),\xi(i)\bigr)=0\qquad \mbox{a.s.},
\]
the family $(\mathcal{R}(k),k\geq1)$ also satisfies this property
with probability $1$, even conditionally given $\mathbf{T}$ and the
Poisson cuts $(t_i,x_i)_{i\in I}$. By Theorem 3 in~\cite{ACRTIII},
this shows that $(\mathcal{R}(k),k\geq1)$ admits a representation as
a continuum random tree, that we call $\operatorname{cut}(\mathbf{T})$. This
means that, still given $\mathbf{T}$ and the process of Poisson marks,
$\operatorname{cut}(\mathbf{T})$ is a pointed metric measure space, with
underlying distance function $r$, root $x_0$ and probability measure
$m$, and that if $x_1,x_2,\ldots$ is an infinite i.i.d. sequence with
distribution $m$, then the matrix $(r(x_i,x_j)\dvtx i,j\geq0)$ has the same
distribution as $(\delta(\xi_i,\xi_j)\dvtx i,j\geq0)$. Up to performing
this ``resampling,'' we thus see that $(\delta(\xi_i,\xi_j)\dvtx i,j\geq
0)$ can itself be seen as the matrix of mutual distances between the
points of an i.i.d. sample of $\operatorname{cut}(\mathbf{T})$. In the
previous discussion, we insisted on conditioning first on $\mathbf{T}$
and $(t_i,x_i)_{i\in I}$ so as to underline the fact that the random
elements $\mathbf{T}$, $(t_i,x_i)_{i\in I}$ and
$\operatorname{cut}(\mathbf{T})$ are defined on a common probability
space. Of course, the equality in distribution (\ref{eqa1}) entails
that unconditionally, the random variable $\operatorname{cut}(\mathbf{T})$
has the same distribution as $\mathbf{T}$.



Let us make a remark at this point. The reader might consider it more
natural to define $\operatorname{cut}(\mathbf{T})$ in a more ``concrete''
way, by first taking the quotient space of $(\mathbf{T}^0,\delta)$ by
the relation $\{(x,y)\dvtx\delta(x,y)=0\}$, and then taking a metric
completion. This operation comes along with a natural mapping that
first projects from $\mathbf{T}^0$ onto the quotient space, and then
injects in the completion. Therefore, we could endow the space with
the image measure of $\mu$ by this natural mapping. However, several
measurability issues appear here: first, this mapping should be
measurable in order that this construction makes sense, and second, it
should be checked that the law of the resulting random metric measure
space is indeed a measurable function of the tree~$\mathbf{T}$. These
issues do not appear in our context because Aldous' construction
defines $\operatorname{cut}(\mathbf{T})$ only in terms of the sequence of
random variables $(\delta(\xi(i),\xi(j))\dvtx i,j\geq0)$. We believe that
these issues can be overcome, however, we are not going to consider
them here to keep this work to a reasonable size.



\subsection{Main result}\label{sec15}

If
$\mathrm{X}=(X,d,\mu,x)$ is a pointed metric measure space and $a>
0$, we let
$a\mathrm{X}=(X,ad,\mu,x)$ be the same space with distances rescaled by
the factor $a$.

%
%
\begin{theorem}\label{T}
\label{sectree-representation}
As $n\to\infty$, we have the following joint convergence in
distribution in $\mathbb{M}\times\mathbb{M}$,\vadjust{\goodbreak} endowed with the
product topology ($\mathbb{M}$ having the Gromov--Prokhorov topology):
\[
\biggl(\frac{\sigma}{\sqrt{n}}{\t}_n,\frac{1}{\sigma\sqrt
{n}}\operatorname{cut}({
\t}_n) \biggr) \Longrightarrow\bigl(\mathbf{T},\operatorname{cut}(
\mathbf{T})\bigr).
\]
Moreover,
$\operatorname{cut}(\mathbf{T})$ has the same distribution as $\mathbf{T}$.
\end{theorem}

It may also be convenient, for example, for readers who would not feel
at ease with weak convergence in the sense induced by the
Gromov--Prokhorov topology, to rephrase the first part of Theorem \ref
{T} as follows (and this is actually what we shall prove).
For every $n\in\N$, set $\xi_n(0)=0$ and consider a sequence $(\xi
_n(i))_{i\geq1}$ of i.i.d. variables having the uniform distribution
on $[n]$. Also let $\xi(0)=0$ and given a CRT ${\mathbf T}$,
let $(\xi(i))_{i\geq1}$ be a sequence of i.i.d. variables in
${\mathbf T}$ distributed according to the mass measure $\mu$. In the
notation of Sections~\ref{sec13} and~\ref{sec14}, we have the following theorem.

%
%
\begin{theorem}\label{T2}
\label{sectree-representation} As $n\to\infty$, the following two
weak convergences hold jointly in the sense of finite-dimensional distributions:
\[
\biggl( \frac{\sigma}{\sqrt{n}} d_n\bigl(\xi_n(i),
\xi_n(j)\bigr)\dvtx i,j\geq1 \biggr) \Longrightarrow\bigl( d\bigl(
\xi(i),\xi(j)\bigr)\dvtx i,j\geq1 \bigr)
\]
and
\[
\biggl( \frac{1}{\sigma\sqrt{n}}\delta_n\bigl(\xi_n(i),
\xi_n(j)\bigr)\dvtx i,j\geq0 \biggr) \Longrightarrow\bigl( \delta\bigl(
\xi(i),\xi(j)\bigr)\dvtx i,j\geq0 \bigr).
\]

\end{theorem}

We conclude this Introduction with three remarks. First, it follows
from Theorem~\ref{T} that
\[
\frac{1}{\sigma\sqrt{n}}\operatorname{cut}({\t}_n) \Longrightarrow
\mathbf{T}.
\]
In particular, the total length $L_k(\operatorname{cut}({\t}_n))$ of the
cut-tree of $\t_n$ reduced to $k$ leaves chosen uniformly at random,
converges after renormalization by a factor $(\sigma\sqrt{n})^{-1}$
to the total length of a Brownian CRT reduced to $k$ i.i.d. leaves
picked according to the mass measure $\mu$. Since the latter is known
to have the Chi$(2k)$-distribution and
$L_k(\operatorname{cut}({\t}_n))$ only differs from the number $N(\bar\t
_n,k)$ of edge-deletions which are needed to recover $k$ distinguished
edges picked uniformly at random in $\bar\t_n$ by at most $k$ units,
we thus see that Proposition~\ref{P1} follows from Theorem~\ref{T}.

Second, we do not know whether the weak convergence stated in the
theorem holds for
stronger topologies of the Gromov--Hausdorff type, even though this is
indeed the case when we only consider the first component.

Third, in order to ease the presentation, it has been convenient to
work with
trees on a set of labeled vertices, namely, $[n]$. This induces a structure
which is not relevant for our main results,\vadjust{\goodbreak} as these could be stated in terms
of graph-theoretic trees over a finite set of vertices
considered up to graph isomorphisms.
The labeling of the vertices comes naturally when considering Galton--Watson
trees; however, it could be ignored after the object is sampled.

\section{Some preliminary results}\label{sec2}
The proof of Theorem~\ref{T} is rather long and relies on several
intermediate results.

\subsection{A modified distance on cut-trees}\label{sec21}
In this section $n$ is fixed; we consider an arbitrary tree $T$ on
$[n]$ and write as usual $\bar T$ for its planted version with $n$
edges. We shall define a modified distance on its cut-tree, which is
both close to the rescaled initial distance on ${\operatorname{cut}}(T)$
and resembles the distance defined in Section~\ref{sec14} for the Brownian
CRT.

Imagine that we mark each edge $e\in\bar T$ with rate $1/\sqrt n$,
independently of the other edges. In particular, the first mark is
assigned after an exponentially distributed time with mean $1/\sqrt n$,
and the edge which is first marked is independent of that time and has
the uniform distribution $\mu_T$ on the set $\bar E$ of the $n$ edges
of $\bar T$ (recall that the set of edges of the planted tree $\bar T$
can be canonically identified with the set $[n]$ of vertices of $T$, so
that by a slight abuse, $\mu_T$ denotes indistinctively the uniform
probability measure on the set of vertices of $T$ or on the set of
edges of $\bar T$).
If we remove the edge $e$ at the instant when it is marked, then we
obtain a continuous time version of edge-deletion process described in
Section~\ref{sec13}. We denote by $\delta_T$ the cut-distance on the
set of
vertices $[n]^0=[n]\cup\{0\}$ which has been defined in that section.
Recall that in this setting, $0$ should be thought of as the root,
$1,\ldots, n$ as leaves, viewed as the edges of $\bar T$, and then
$\delta_{T}(0,i)$ is given by the number of edge-removals that are
performed on the successive blocks containing $i$ until $i$ is finally removed.

For every
$t\geq0$, we write $\bar T(t)$ for the random forest that results from
the edge deletion process at time $t$, and
for every $i\in[n]$, $\bar T_i(t)$ for the tree-component of
$\bar T(t)$ which contains the edge labeled by $i$, agreeing, of
course, that $\bar T_i(t)$ is empty whenever $i$ has been removed
before time $t$.
We also set $\mu_{T,i}(t)=\mu_T(\bar T_i(t))$; this quantity
gives the number of edges of $\bar T_i(t)$ normalized by a factor $1/n$.
Mimicking the construction of the cut-distance on the Brownian CRT in
Section~\ref{sec14}, we now introduce
\[
\delta'_T(0,0)=0,\qquad \delta'_T(0,
i)= \delta'_T(i,0) = \int_{0}^\infty
\mu_{T,i}(t)\,\d t,\qquad i\in[n],
\]
and
\[
\delta'_T(i,j) = \int_{t_{ij}}^\infty
\bigl(\mu_{T,i}(t)+\mu_{T,j}(t)\bigr)\,\d t,\qquad i,j\in[n],
\]
where $t_{ij}$ is the first instant $t$ at which the edges $i$ and $j$
become disconnected in $\bar T(t)$.

We have thus endowed $[n]^0$ with two distances, $\delta_T$ and
$\delta'_T$ related to the edge-deletion process;
our purpose here is to observe that $\frac{1}{\sqrt n} \delta_T$ and
$\delta'_T$
are close when $n$ is large.
Here is a precise statement.

%
%
\begin{lemma} \label{lem1} For every $i\in[n]$, we have
\[
\E\biggl(\biggl\llvert\frac{1}{\sqrt n} \delta_T(0,i)-
\delta'_T(0,i)\biggr\rrvert^2 \biggr) =
\frac{1}{\sqrt n} \E\bigl(\delta'_T(0,i)\bigr)
\]
and, as a consequence, for every $j\in[n]$, we also have
\[
\E\biggl(\biggl\llvert\frac{1}{\sqrt n} \delta_T(i,j)-
\delta'_T(i,j)\biggr\rrvert^2 \biggr) \leq
\frac{2}{\sqrt n} \E\bigl(\delta'_T(0,i) +
\delta'_T(0,j)\bigr).
\]

\end{lemma}
\begin{pf}
%
We shall focus on the first inequality as the second can be established
using a closely related argument.

Denote by $N^i_{T}(t)$ the number of edges that have been removed up to
time $t$ from the tree-components that contain the edge $i$; in particular,
\[
\lim_{t\to\infty} N^i_{T}(t) = \delta_T(0,i).
\]
Since each edge of $\bar T(t)$ is removed with rate $1/\sqrt n$,
independently of the other edges, the process
\[
M(t)= \frac{1}{\sqrt n} N^i_{T}(t)- \int
_0^t \mu_{T,i}(s) \,\d s ,\qquad t\geq0,
\]
is a purely discontinuous martingale with terminal value
\[
\lim_{t\to\infty}M(t) =\frac{1}{\sqrt n} N^i_{T}(
\infty) - \int_0^{\infty}\mu_{T,i}(s) \,\d s =
\frac{1}{\sqrt n} \delta_T(0,i) - \delta'_T(0,i).
\]
Further, its quadratic variation is $[M]_t=n^{-1} N^i_{T}(t)$
and thus its oblique bracket is given by
\[
\langle M\rangle_t=\frac{1}{\sqrt n} \int_0^t
\mu_{T,i}(s) \,\d s.
\]
As a consequence, we have
\[
\E\biggl( \biggl\llvert\frac{1}{\sqrt n} \delta_T(0,i) -
\delta'_T(0,i) \biggr\rrvert^2 \biggr) =
n^{-1/2} \E\biggl(\int_0^{\infty}
\mu_{T,i}(s) \,\d s \biggr),
\]
which is our statement.
\end{pf}

\subsection{Joint convergence of the subtree sizes}\label{sec22}

Recall that $\t_n$ denotes a Galton--Watson tree corresponding to
critical reproduction law with finite variance $\sigma^2>0$ and
conditioned to have size $n$.
We know from Aldous~\cite{ACRTIII} that\vadjust{\goodbreak} $\sigma n^{-1/2} \t_n$
converges in distribution to a Brownian CRT.
Motivated by Lemma~\ref{lem1}, the purpose of this section is to point
out that this convergence can be reinforced to hold jointly with that
of the rescaled sizes of subtrees appearing in the edge-deletion processes.

Before providing a rigorous statement, we need to introduce some
further notation.
For $n$ fixed, given the random tree $\t_n$, we consider a sequence
$(\xi_n(i), i\geq1)$ of i.i.d. random variables in $[n]$, each having
the uniform distribution $\mu_n$ on $[n]$. We stress that the $\xi
_n(i)$ should be viewed as random edges of the planted tree
$\bar\t_n$, although in this section it will sometimes be convenient
to think of the latter as vertices of $\t_n$. Randomly marking each
edge of $\bar{\t}_n$ at rate $1/\sqrt n$ as in Section~\ref{sec21}, we denote
for every $t\geq0$ by $\mu_{n,\xi_n(i)}(t)$ the
number of edges of the tree-component containing the edge $\xi_n(i)$
in the forest at time $t$,
$\bar{\t}_n(t)$, and normalized by a factor $1/n$.
We also denote by $\tau_n(i,j)$ the first instant when the edges $\xi
_n(i)$ and $\xi_n(j)$ are disconnected in the forest $\bar{\t}_n(t)$.

Next, we consider the Brownian CRT ${\mathbf T}$ together with the
Poisson point process of marks on its skeleton as in Section \ref
{sec14}, and
an independent sequence
$(\xi(i), i\geq1)$ of i.i.d. random variables in ${\mathbf T}$
distributed according to the uniform measure~$\mu$. Recall that for
every $t\geq0$, $\mu_{\xi(i)}(t)$ denotes the $\mu$-mass of the
tree-component containing $\xi(i)$ in forest
${\mathbf T}(t)$ that results from ${\mathbf T}$ by cutting its
skeleton at marks which appeared before time $t$. Finally, we denote by
$\tau(i,j)=t_{\xi(i)\xi(j)}$ the first
instant when the points $\xi(i)$ and $\xi(j)$ are disconnected, that
is, the first time when a mark is put on the segment joining $\xi(i)$
and $\xi(j)$ in ${\mathbf T}$.

We are now able to state the following lemma.

%
%
\begin{lemma}\label{lem2} As $n\to\infty$, we have the following
weak convergences:
\begin{eqnarray*}
\frac{\sigma}{\sqrt{n}}{\t}_n &\Longrightarrow&\mathbf{T},
\\
\bigl(\tau_n(i,j)\dvtx i,j\in\N\bigr) &\Longrightarrow&\bigl(\tau
(i,j)\dvtx i,j\in\N\bigr)
\end{eqnarray*}
and
\[
\bigl( \mu_{n,\xi_n(i)}(t)\dvtx t\geq0 \mbox{ and }i\in\N\bigr)
\Longrightarrow
\bigl( \mu_{\xi(i)}(t)\dvtx t\geq0 \mbox{ and }i\in\N\bigr),
\]
where the three hold jointly; the first in the sense induced by the
Gromov--Prokhorov topology, and the second and third in the sense of
finite-dimensional distributions.

\end{lemma}
\begin{pf}
The proof closely follows arguments developed by Aldous and Pitman
\cite{AP} in a similar setting (see Section 2.3 there).

We first consider the edge-deletion process on $\t_n$ rather than on
its planted version $\bar\t_n$ and
view the random variables $(\xi_n(i), i\geq1)$ as a sequence of
vertices of $\t_n$
rather than edges of $\bar{\t}_n$. Let ${\mathcal R}(n,k)$ denote the
subtree of $\t_n$ spanned by the first $k$ random vertices $\{\xi_n(i),
1\leq i \leq k\}$ and the root of $\t_n$. Similarly, we denote
by ${\mathcal R}(\infty, k)$ the subtree obtained from the CRT
${\mathbf T}$ by reduction to its root and the first $k$ i.i.d.
variables $\{\xi(i), 1\leq i \leq k\}$ with common distribution
the mass-measure $\mu$ on ${\mathbf T}$.
Here, we adopt the framework of Aldous~\cite{ACRTIII}, viewing the
reduced trees as a combinatorial rooted tree structure with edge
lengths and leaves labeled by $1,\ldots, k$.
As it was already stressed, we know from the work of Aldous \cite
{ACRTIII} that there is the convergence
%
%
\begin{equation}
\label{eqct} \frac{\sigma}{\sqrt{n}}{\t}_n \Longrightarrow\mathbf{T},
\end{equation}
and this can be rephrased in terms of reduced trees as
%
%
\begin{equation}
\label{eqcrt} \frac{\sigma}{\sqrt{n}}{\mathcal R}(n,k) \Longrightarrow
{\mathcal R}(
\infty,k),
\end{equation}
where $\frac{\sigma}{\sqrt{n}}{\mathcal R}(n,k)$ has the same
tree-structure as ${\mathcal R}(n,k)$ but with edge lengths rescaled by
a factor $\sigma/\sqrt{n}$. We shall now see how (\ref{eqct}) can be
enriched to encompass the further convergences in the statement.

Next, we write $({\mathcal R}(n,k,t)\dvtx t\geq0)$ for the reduced tree
${\mathcal R}(n,k)$ endowed with a point process of marks on its edges.
More precisely, each edge receives a mark at its midpoint precisely at
the time when this edge is removed as in Section~\ref{sec21}.
Similarly, we
denote by
$({\mathcal R}(\infty,k,t)\dvtx t\geq0)$ for the reduced tree ${\mathcal
R}(\infty,k)$ endowed with a Poisson point process of marks on its
skeleton with intensity $\d t \otimes\d\lambda$, where by slightly
abusive notation, $\lambda$ is now the length measure on the reduced
tree ${\mathcal R}(\infty,k)$. It is then easy to extend (\ref
{eqcrt}) to
%
%
\begin{equation}
\label{eqcrtm} \biggl(\frac{\sigma}{\sqrt{n}}{\mathcal R}(n,k, t)\dvtx
t\geq
0 \biggr)
\Longrightarrow\bigl({\mathcal R}(\infty,k, t/\sigma)\dvtx t\geq0 \bigr)
\end{equation}
on the space of rooted trees with $k$ leaves and edge-lengths, endowed
with an increasing process of marked points, this space being equipped
with the appropriate topology.
More precisely, the time-rescaling with a factor $1/\sigma$ in the
right-hand side stems from the fact that the edges in ${\mathcal
R}(n,k)$ have been rescaled by $\sigma/\sqrt n$.

Then, for every $i\geq1$, denote by $\eta(n,k,i,t)$ the number of
vertices among $\xi_n(1), \ldots, \xi_n(k)$ in the tree-component
containing the vertex $\xi_n(i)$ which results from ${\mathcal R}(n,k)$
by cutting at marks that appeared before time $t$.
Denote also by $\tau'_n(i,j)$ the first instant when a mark appears on
the segment in ${\mathcal R}(n,k)$ connecting $\xi_n(i)$ and $\xi_n(j)$.

Similarly, let $\eta(\infty,k,i,t)$ be the number of vertices among
$\xi(1), \ldots, \xi(k)$ in the tree-component containing the vertex
$\xi(i)$ which results from ${\mathcal R}(\infty,k)$ by cutting at
marks that
appeared before time $t$. It follows from (\ref{eqcrtm}) that
\[
\bigl(\eta(n,k,i,t)\dvtx t\geq0 \mbox{ and }i\in\N\bigr) \Longrightarrow
\bigl(
\eta(\infty,k,i,t/\sigma)\dvtx t\geq0 \mbox{ and }i\in\N\bigr)
\]
and
%
%
\begin{equation}
\label{eqtemps} \bigl(\tau'_n(i,j)\dvtx i,j\in\N\bigr)
\Longrightarrow\bigl(\sigma\tau(i,j)\dvtx i,j\in\N\bigr)
\end{equation}
in the sense of finite-dimensional distributions. More precisely, these
convergences hold jointly,
also together with (\ref{eqct}).

The law of large numbers entails that for each fixed $i$ and $t\geq0$,
\[
\lim_{k\to\infty} k^{-1} \eta(\infty,k,i,t/\sigma) =
\mu_{\xi
(i)}(t/\sigma)\qquad \mbox{almost surely.}
\]
We deduce that for every fixed integer $\ell$ and times $0\leq t_1
\leq\cdots\leq t_{\ell}$, we can construct a sequence $k_n\to\infty
$ sufficiently slowly, such that
\[
\bigl(k_n^{-1}\eta(n,k_n,i,t_j)\dvtx
1\leq i, j \leq\ell\bigr) \Longrightarrow\bigl(\mu_{\xi(i)}(t_j/
\sigma)\dvtx1\leq i, j \leq\ell\bigr),
\]
and again this convergence holds jointly with (\ref{eqct}) and (\ref
{eqtemps}).

This is essentially the sought-after result; the only minor difference
is that we viewed the $\xi_n(i)$ as vertices of $\t_n$ rather than
edges of the planted tree $\bar{\t}_n$. However, it is easy to check
that, with high probability, this makes no difference when $n$ is
large. Indeed, we realize that on the event that the edges $\xi_n(i)$
and $\xi_n(j)$ of $\bar{\t}_n$ have not been marked before time $t$,
which has probability greater than $2\exp(-t/\sqrt n)-1\to1$ as $n\to
\infty$,
$ \eta(n,k,i,t) $ differs from $k\mu_{n,\xi_n(i)}(t)$ by at most one
unit (recall that a tree with $j$ vertices has $j-1$ edges).
Similarly, on the event that the edges $\xi_n(i)$ and $\xi_n(j)$ have
not been removed when the segment connecting $\xi_n(i)$ and $\xi_n(j)$
receives its first mark, which has probability
$\E(d_n(\xi_n(i), \xi_n(j))/(2+d_n(\xi_n(i), \xi_n(j))
) \to1$ as $n\to\infty$,
there is the identity $\tau_n(i,j)=\tau'_n(i,j)$.
The proof is complete.
\end{pf}

\subsection{A uniform bound}\label{sec23}
The next technical step in the proof of Theorem~\ref{T} is the
obtention of a uniform bound
for the first moment of the size of a ``typical'' tree component
occurring in the
edge-deletion process for Galton--Watson trees.
Specifically, recall that $\xi_n=\xi_n(1)$ is a random edge of $\bar
{\t}_n$ picked according to the uniform probability measure $\mu_n$
and independently of the edge-deletion process, and that
$\mu_{n,\xi_n}(t)$ denotes the number of edges in the tree component
of $\bar{\t}_n(t)$ which contains the random edge $\xi_n$ and
rescaled by a factor $1/n$.
We claim the following.

%
%
\begin{lemma} \label{lem3}
There exists some finite constant $C>0$ depending only on the offspring
distribution $\nu$
such that
\[
\E\bigl(\mu_{n,\xi_n}(t)\bigr) \leq C\frac{\exp(-t/\sqrt n)}{n(1-\exp
(-t/\sqrt n))^2}
\]
for all $t\geq0$ and $n\in\N$.
\end{lemma}
We stress that this bound is only relevant when $t$ is not too small,
since the left-hand side is always less than or equal to $1$.\vadjust{\goodbreak}

The proof of Lemma~\ref{lem3} relies crucially on an estimate due to
Janson~\cite{Janson} and an invariance property under random
re-rooting for planted Galton--Watson trees. It is convenient to
postpone its proof to Section~\ref{sec4}; we merely conclude this
section with a
consequence of that lemma which will be used in the proof of Theorem
\ref{T}.

%
%
\begin{corollary}\label{C1} There exists a finite constant $C'>0$
depending only on the offspring distribution $\nu$ such that
\[
\E\bigl(\delta'_n(\xi_n,0)\bigr) \leq
C' \qquad\mbox{for all }n\in\N,
\]
where $\delta'_n$ denotes the modified distance on
$\operatorname{cut}(\t_n)$ defined in Section~\ref{sec21}.
\end{corollary}

\begin{pf}
An application of Lemma~\ref{lem3} at the second line below gives
\begin{eqnarray*}
\E\bigl(\delta'_n(\xi_n,0)\bigr) &=& \int
_0^{\infty}\E\bigl(\mu_{n,\xi
_n}(t) \bigr) \,\d t
\\
&\leq& 1+ C \int_1^{\infty} \frac{\exp(-t/\sqrt n)}{n(1-\exp
(-t/\sqrt n))^2} \,\d t
\\
&=& 1+ \frac{C} {{\sqrt n}(1-\exp(-1/\sqrt n))},
\end{eqnarray*}
and this last quantity remains indeed bounded as
$n\to\infty$.
\end{pf}

\section{\texorpdfstring{Proof of Theorem \protect\ref{T}}{Proof of Theorem 1}}\label{sec3}

\subsection{The cut-tree of a Brownian CRT is another Brownian
CRT}\label{sec31}\label{seccut-tree-brownian}
In this section, we complete the construction of
$\operatorname{cut}(\mathbf{T})$ that was performed in Section
\ref{seccut-tree-continuum}, and show that it has the same
distribution as the Brownian CRT. Both will follow from the following
lemma.

%
%
\begin{lemma}\label{lem5} Set
$\xi(0) \equiv0$ and let $(\xi(i)\dvtx i\in\N)$ denote a sequence of
i.i.d. points in $\mathbf{T}$ distributed according to the uniform
probability measure $\mu$. Then there is the identity in law
%
%
\begin{equation}
\label{eqid} \bigl(d\bigl(\xi(i+1),\xi(j+1)\bigr)\dvtx i,j\geq0 \bigr)
\stackrel{(\mathit{law})}{=} \bigl(\delta\bigl(\xi(i),\xi(j)\bigr)\dvtx
i,j\geq0
\bigr).
\end{equation}
\end{lemma}

As a warmup, we first provide a short proof of the one-dimensional
identity in~(\ref{eqid}), that is, for $i$ and $j$ fixed.
By a well-known property of invariance in law of the Brownian CRT under
random uniform re-rooting, it suffices to treat the case $i=0$ and $j=1$,
and we thus consider the cut-distance $\delta(0,\xi)$ of a random
point $\xi$ with law $\mu$
to the root $\xi(0)=0$. In the notation of Section~\ref{sec14}, we have
\[
\delta(0,\xi)= h_{\xi}=\int_0^{\infty}
\mu_{\xi}(t) \,\d t.
\]
Applying the identity in distribution (\ref{eq1}), we see that this
variable has the same law as
the Cauchy transform
\[
C(\sigma)=\int_0^{\infty} \frac{\d t}{1+\sigma(t)},
\]
where $(\sigma(t), t\geq0)$ is the stable$(1/2)$ subordinator given
by the first-passage time process of a standard Brownian motion, that
is, with Laplace exponent $\Phi(r)=\sqrt{2r}$.

According to Corollary 3 and Lemma 3 of~\cite{Be2} specified to
subordinators, we have
\[
\P\bigl(C(\sigma)\leq t\bigr)=1-\exp\bigl(-\gamma(t)\bigr),\qquad t\geq0,
\]
where $\gamma$ denotes the inverse of the function $t\to\int_0^t
\Phi(r)^{-1}\,\d r$.
For $\Phi(r)=\sqrt{2r}$, we get $\gamma(t)=\frac{1}{2}t^2$ and
conclude that the distribution function of
$C(\sigma)$ is $t\to1-\exp(-\frac{1}{2}t^2)$, which is the
distribution function of the Rayleigh law.
The latter coincides with the distribution of the height $d(0,\xi)$ of
a point picked uniformly at random in $\mathbf{T}$, and we conclude
that indeed (\ref{eqid}) holds in the weaker sense of one-dimensional
distributions. We note passing by that the claim (\ref{eq2}) is now
established.

Unfortunately, such direct calculations are not available for
multidimensional distributions, and we shall use a different approach
which relies on a general feature of self-similar fragmentations.
We thus start by developing elements in this area and refer the reader
to~\cite{RFC} and, in particular, Chapters~2 and 3 there for background.

For this purpose, it is convenient to work in the setting of
processes with values in the space of partitions of $\N=\{1,2,
\ldots\}$, which arise naturally from i.i.d. sampling. Given
$\mathbf{T}$ and a sequence $(\xi(i)\dvtx i\in\N)$ of i.i.d. points with
law $\mu$, we shall consider two such fragmentation processes. A
first fragmentation process $\Gamma=(\Gamma(t), t\geq0)$ results
from cutting the CRT at its heights. Specifically, recall that $0$
denotes the root of $\mathbf{T}$. For every $x,y\in\mathbf{T}$, let
$[x,y]$ be the segment connecting $x$ and $y$, and define the
branch-point $x\wedge y$ as the unique point in $\mathbf{T}$ such
that $[0,x]\cap[0,y]=[0, x\wedge y]$. Then we declare that two
distinct integers $i\neq j$ belong to the same block of the partition
$\Gamma(t)$ if and only if the height of the branch-point of $\xi(i)$
and $\xi(j)$ is greater than $t$, that is, $d(0,\xi(i)\wedge
\xi(j))>t$. In other words, the height of the branch-point is the
first time at which $i$ and $j$ are disconnected in the fragmentation
process~$\Gamma$. We stress that $\{i\}$ is a singleton of
$\Gamma(t)$ whenever the height of $\xi(i)$ is smaller than or equal
to $t$; in particular, $\Gamma(t)$ is the partition into singletons
whenever $t \geq\sup\{d(0,x)\dvtx x\in\mathbf{T}\}$.

Recall that $(t_i, x_i)_{i\in I}$ denotes the family of the atoms of a
Poisson random
measure with intensity $\d t\otimes\d\lambda$ on the skeleton of
$\mathbf{T}$, which is assumed to be independent of the preceding.
We denote by $\Pi=(\Pi(t), t\geq0)$ the Aldous--Pitman fragmentation
of the Brownian CRT,
obtained by declaring that two integers $i$ and $j$ belong to the same
block of $\Pi(t)$ if and only if
$[\xi(i), \xi(j)] \cap\{x_i\dvtx t_i\leq t\}=\varnothing$, that is, if
and only if $\xi(i)$ and $\xi(j)$ belong to the same component of the
random forest $\mathbf{T}(t)$.

These two fragmentation processes are related by a sort of time-substitution
which is the key to Lemma~\ref{lem5}. For every $i\in\N$ and $t\geq
0$, denote by $B_i(t)$ the block of $\Pi(t)$ which contains $i$
and by $|B_i(t)|$ its asymptotic frequency; it is also convenient to
agree that $B_i(\infty)=\{i\}$.
Next define
\[
\rho_i(t)=\inf\biggl\{u\geq0\dvtx\int_0^u
\bigl|B_i(r)\bigr| \,\d r >t \biggr\},\qquad t\geq0,
\]
with the usual convention $\inf\varnothing= \infty$. Roughly
speaking, we use the $\rho_i$ to time-change the fragmentation $\Pi$
and write $\Pi'(t)$ for the partition whose family of blocks is given
by the $B_i(\rho_i(t))$ for $i\in\N$ (observe that two such blocks are
either disjoint or equal).

%
%
\begin{lemma}\label{lem6} In the notation above, the fragmentation processes
$\Gamma$ and $\Pi'$ have the same law.
\end{lemma}
\begin{pf}
The Aldous--Pitman fragmentation $\Pi$ is a self-similar fragmentation
with index $\alpha=1/2$, erosion coefficient $0$ and dislocation
measure denoted here by~$\Delta$,
as it is seen from, for example, Theorem 3 in~\cite{AP} and Theorem
5.4 in~\cite{RFC}.
According to Theorem 3.3 in~\cite{RFC}, the time-changed fragmentation
$\Pi'=(\Pi(\rho(t)), t\geq0)$
is then a self-similar fragmentation, now with index $\alpha-1=-1/2$,
with no erosion and the same dislocation measure $\Delta$.

On the other hand, the discussion in~\cite{Be1}, pages 339 and 340, and
the well-known
construction of the Brownian CRT from twice the normalized Brownian
excursion (see, e.g., Corollary 22 in~\cite{ACRTIII}) show
that $\Gamma$ is again a self-similar fragmentation
with index $-1/2$, no erosion and dislocation measure $\Delta$.
Hence, $\Gamma$ and $\Pi'$ are two self-similar fragmentations with
the same characteristics; they thus have the same law (see~\cite{RFC},
page 150).
\end{pf}

Lemma~\ref{lem5} should now be obvious. Indeed, by the law of large
numbers, the $\mu$-mass of a component of $\mathbf{T}(t)$ can be
recovered as the asymptotic frequency of the corresponding block of the
partition, and, in particular, $|B_i(t)|=\mu_{\xi(i)}(t)$.
Recall that the height $d(0,\xi(i))$ of $\xi(i)$ in $\mathbf{T}$ can
be viewed as the first instant $t$ when
$\{i\}$ is a singleton of $\Gamma(t)$, a quantity which, in terms of
the Aldous--Pitman fragmentation~$\Pi$, corresponds to
\[
\int_0^{\infty}\bigl|B_i(t)\bigr| \,\d t = \int
_0^{\infty}\mu_{\xi(i)}(t) \,\d t = \delta\bigl(0,
\xi(i)\bigr).
\]
Similarly, for $i,j\in\N$ with $i\neq j$,
\[
d\bigl(0,\xi(i)\wedge\xi(j)\bigr) = \tfrac{1}{2} \bigl(d\bigl(0, \xi(i)
\bigr) + d\bigl(0, \xi(j)\bigr) - d\bigl( \xi(i),\xi(j)\bigr) \bigr),
\]
and in terms of $\Pi$, the last quantity corresponds to
\[
\int_0^{\tau(i,j)}\bigl|B_i(t)\bigr| \,\d t = \int
_0^{\tau(i,j)}\mu_{\xi
(i)}(t) \,\d t = \delta\bigl(0,
\xi(i)\wedge\xi(j)\bigr),
\]
where $\tau(i,j)$ denotes the first instant $t$ when a mark appears on
the segment $[\xi(i), \xi(j)]$.
Combining these observations with Lemma~\ref{lem6}, we conclude that
(\ref{eqid}) holds.


\subsection{Proof of weak convergence}\label{sec32}
It is convenient to first establish the convergence in Theorem~\ref{T}
when $\operatorname{cut}(\t_n)$ is endowed with the modified distance
$\delta'_n$ as defined in Section~\ref{sec21}. We write
$\operatorname{cut}'(\t_n)=([n]^0, \delta'_n, \mu_n, 0)$ for the
pointed metric measure
space equipped with the modified distance and claim the following.

%
%
\begin{lemma} \label{lem4}
As $n\to\infty$, there is the joint convergence in
the weak sense induced by the Gromov--Prokhorov topology
\[
\biggl(\frac{\sigma}{\sqrt{n}}{\t}_n,\operatorname{cut}'({
\t}_n) \biggr) \Longrightarrow\bigl(\mathbf{T}, \sigma\operatorname{cut}(
\mathbf{T})\bigr).
\]
\end{lemma}
\begin{pf}
We use the setting and notation of Section~\ref{sec22} and derive from Lemma
\ref{lem2} that
for every fixed integer $\ell$,
\[
\frac{\sigma}{\sqrt{n}}{\t}_n \Longrightarrow\mathbf{T}
\]
and
\[
\Biggl( 2^{-\ell} \sum_{j=1}^{4^{\ell}}
\mu_{n,\xi_n(i)}\bigl(j 2^{-\ell
}\bigr)\dvtx i\in\N\Biggr) \Longrightarrow
\Biggl( 2^{-\ell} \sum_{j=1}^{4^{\ell}}
\mu_{\xi(i)}\bigl(j 2^{-\ell
}/\sigma\bigr)\dvtx i\in\N\Biggr),
\]
where the two convergences hold jointly; the first in the sense induced
by the Gromov--Prokhorov topology and the second in the sense of
finite-dimensional distributions.

For every nonincreasing function $f\dvtx[0,\infty)\to[0,1]$, there is
the bound
\[
\Biggl\llvert\sigma\int_0^{\infty} f(t)\,\d t -
2^{-\ell} \sum_{j=1}^{4^{\ell}}f\bigl(j
2^{-\ell}/\sigma\bigr)\Biggr\rrvert\leq\sigma\biggl(2^{-\ell} +
\int_{2^{\ell}/\sigma}^{\infty} f(t)\,\d t \biggr).
\]
Since the Rayleigh distribution has a finite mean, we deduce from
(\ref{eq2}) that
\[
\E\Biggl(\Biggl\llvert\sigma\int_0^{\infty}
\mu_{\xi(i)}(t)\,\d t - 2^{-\ell} \sum_{j=1}^{4^{\ell}}
\mu_{\xi(i)}\bigl(j 2^{-\ell}/\sigma\bigr)\Biggr\rrvert\Biggr) \to
0 \qquad\mbox{as } \ell\to\infty,
\]
where, of course, the left-hand side does not depend on $i$.

Similarly, now using Lemma~\ref{lem3}, we obtain the uniform bound
\begin{eqnarray*}
&&
\E\Biggl(\Biggl\llvert\int_0^{\infty}
\mu_{n,\xi_n(i)}(t)\,\d t - 2^{-\ell} \sum_{j=1}^{4^{\ell}}
\mu_{n,\xi_n(i)}\bigl(j 2^{-\ell
}\bigr)\Biggr\rrvert\Biggr) \\
&&\qquad \leq
2^{-\ell} + C\int_{2^{\ell}}^{\infty}
\frac{\exp
(-t/\sqrt n)}{n(1-\exp(-t/\sqrt n))^2} \,\d t
\\
&&\qquad = 2^{-\ell} + \frac{C} {{\sqrt n}(1-\exp(-2^{\ell}/\sqrt n))}.
\end{eqnarray*}
We thus see that the left-hand side above also tends to $0$ as $\ell
\to\infty$, uniformly in $n\in\N$ (and again these quantities do
not depend on $i$).

Recalling that
\[
\delta'_n\bigl(0,\xi_n(i)\bigr)= \int
_0^{\infty} \mu_{n,\xi_n(i)}(t)\,\d t \quad\mbox{and}\quad
\delta\bigl(0,\xi(i)\bigr)= \int_0^{\infty}
\mu_{\xi(i)}(t)\,\d t,
\]
we conclude that
\[
\bigl( \delta'_n\bigl(0,\xi_n(i)\bigr)\dvtx i
\in\N\bigr) \Longrightarrow\bigl( \sigma\delta\bigl(0,\xi(i)\bigr)\dvtx
i\in\N
\bigr)
\]
in the sense of finite-dimensional distributions, and the latter holds
jointly with
$\frac{\sigma}{\sqrt{n}}{\t}_n
\Longrightarrow\mathbf{T}$.
Essentially, the same argument, now further using the convergence of
disconnection times stated in Lemma~\ref{lem2}, shows that the
preceding also hold jointly with
\[
\bigl( \delta'_n\bigl(\xi_n(i),
\xi_n(j)\bigr)\dvtx i,j\in\N\bigr) \Longrightarrow\bigl( \sigma\delta
\bigl(\xi(i),\xi(j)\bigr)\dvtx i,j\in\N\bigr).
\]
This is
precisely the meaning of our statement, since we have seen in
Section~\ref{seccut-tree-continuum} that the doubly-infinite sequence
$(\delta(\xi(i),\xi(j))\dvtx i,j\geq0)$ can be seen as the matrix of
mutual distances between the root of $\operatorname{cut}(\mathbf{T})$
and an
i.i.d. sample of points in $\operatorname{cut}(\mathbf{T})$.
\end{pf}

This immediately entails the convergence stated in Theorem~\ref{T}.
Specifically, recall that $\xi_n(0)\equiv0$; combining Lemma \ref
{lem1} and Corollary~\ref{C1}, we get that for all $i,j\geq0$ there
is the upper-bound
\[
\E\biggl(\biggl\llvert\frac{1}{\sqrt n} \delta_n\bigl(
\xi_n(i), \xi_n(j)\bigr)- \delta'_n
\bigl(\xi_n(i), \xi_n(j)\bigr)\biggr
\rrvert^2 \biggr) \leq\frac
{4C'}{\sqrt n}.
\]
Therefore, Lemma~\ref{lem4} can be rephrased as
\[
\biggl(\frac{\sigma}{\sqrt{n}}{\t}_n, \frac{1}{\sqrt{n}}\operatorname{cut}({
\t}_n) \biggr) \Longrightarrow\bigl(\mathbf{T}, \sigma\operatorname{cut}(
\mathbf{T})\bigr),
\]
which is the claimed convergence.

\section{\texorpdfstring{Proof of Lemma \protect\ref{lem3}}{Proof of Lemma 3}}\label{sec4}
The purpose of this final section is to establish Lem\-ma~\ref{lem3}.
The proof relies on an estimate due to Janson~\cite{Janson}, combined
with an invariance property of the law of Galton--Watson\vadjust{\goodbreak} trees under
random re-planting. We start by the latter; this is the main part of
this paper where working with planted trees makes the approach simpler.

It will be convenient in this section to work in the setting of planar
rooted trees
rather than tree structures on a set of labeled vertices. As vertices
of planar rooted trees
can be canonically enumerated, for instance, in the breadth-first
search order,
the transformations appearing in this section could also be re-phrased
in terms of tree structures on
a set of labeled vertices, though doing so would make the descriptions
more involved.

A Galton--Watson tree $\t$ is thus viewed here as a random planar
rooted tree. We write $V(\t)$ for the set of vertices of $\t$; in
particular, its cardinal
$|V(\t)|$ is the total number of individuals in the branching process
with critical reproduction law $\nu$ and whose genealogy is
represented by $\t$.
Recall also that $\bar\t$ denotes the planted version of $\t$ and
then $|V(\t)|$ is the number of edges of $\bar\t$.\vspace*{1pt}

A \textit{pointed tree} is a pair $(\bar T,v)$ where $\bar T$ is a planted
planar tree and $v$ a vertex distinct from the base, that is, a vertex
of $T$.
We endow the space of pointed trees with a sigma-finite measure
$\mathrm{GW}_*$ defined by
\[
\mathrm{GW}_*(\bar T,v)= \P(\bar\t= \bar T),
\]
where $\bar T$ denotes a generic planar planted tree and $v\in V(T)$.
This measure is a classical object appearing, in particular, in the
approach by Lyons, Pemantle and Peres~\cite{LPP}.

\begin{figure}

\includegraphics{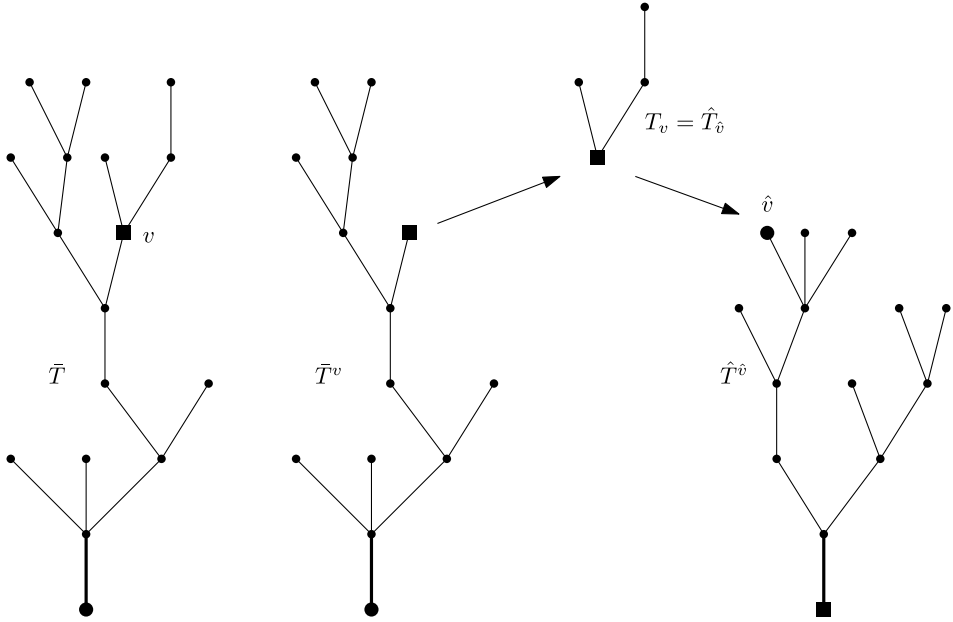}

\caption{The trees $\bar{T}_v,\bar{T}^{v},\hat{T}^{\hat v}$, with
the distinguished vertices $v$ and $\hat{v}$, with the last picture
explaining how to construct $\hat{T}$ from $\hat{T}^{\hat v}$ and $T_v$.}
\label{figrerooting}
\end{figure}

We now describe a transformation of pointed trees which will be used
in the proof of Lemma~\ref{lem3}. If $(\bar T,v)$ is a pointed tree,
we let $T_v$ be the (nonplanted) tree formed by all the descendants
of $v$ in $T$, including $v$, and $\bar T^v$ be the subtree obtained
by removing all the strict descendants of $v$ in $\bar T$. We first
re-plant $\bar T^v$ at $v$, viewing the edge connecting $v$ to its
parent in $\bar T^v$ as the new base-edge. We denote the former
base-vertex by $\hat v$ and the new planted tree by $\hat{T}^{\hat
v}$. Finally, we re-graft $T_v$ at $\hat v$ and get another pointed
tree which we denote by $(\hat T, \hat v)$ (see Figure
\ref{figrerooting}).

%
%

We stress that the set of vertices of $\bar T$ and of $\hat T$
coincide. More precisely, if we remove the strict descendants of $\hat
v$ in $\hat T$, then we get $\hat T^{\hat v}$
(so our notation is coherent), while the subtree formed by the
descendants of $\hat v$ in $\hat T$ coincides with $T_v$, that is,
$\hat T_{\hat v}= T_v$. We also\vspace*{1pt} note that the transformation of pointed trees
$(\bar T,v)\to(\hat T, \hat v)$ is involutive, that is, its iteration
is the identity.

%
%
\begin{proposition}
\label{secverif-hypoth-eqrefhy}
The ``laws'' of $(\hat{T},\hat{v})$
and $(T,v)$ under the measure
$\mathrm{GW}_*$ are the same. Equivalently,
$(\hat{T}^{\hat v},T_v)$ and $(T^{v},T_v)$ have same ``law''
under $\mathrm{GW}_*$.
\end{proposition}
\begin{pf}
We only sketch the proof of this proposition, leaving some technical
details to the
reader. From Lyons, Pemantle and Peres~\cite{LPP}, we know that under
$\mathrm{GW}_*$, a typical pointed tree\vadjust{\goodbreak} $(\bar T,v)$ can be described
in terms of Kesten's critical
Galton--Watson tree conditioned on nonextinction, or size-biased
tree. We start by recalling some features concerning the latter.

The size-biased tree is a random planted tree $\t_*$ with a single infinite
branch $v_0,v_1,\ldots$ starting from the base $v_0$, that can be
constructed as follows. As for a standard Galton--Watson tree, every
vertex in $\t_*$ has an offspring that
is distributed according to the reproduction law $\nu$ and
independently of the other vertices, except, of course, the base $v_0$
which has only one child, and the further vertices
of the infinite branch, $v_1, v_2, \ldots\,$, whose offspring
distribution is the
size-biased measure $\nu_*(k)=k\nu(k)$ (recall that the reproduction
law $\nu$ is critical). The size-biased tree $\t_*$ is constructed
inductively starting from the base $v_0$ by claiming that the
$i$th vertex $v_i$ on the infinite branch is chosen uniformly at random
from the
offspring of $v_{i-1}$.

For every $h\geq1$, let $\t_*^{v_h}$ be the planted tree
pointed at $v_h$ which is obtained from $\t_*$ by removing all the strict
descendants of $v_h$, and let $\mathrm{GW}_*^h$ be its law. Now
consider an integer $h\geq1$
``sampled'' according to the counting
measure on~$\N$. Sample the pointed tree $(\t_*^{v_h}, v_h)$ as
above, and
independently choose a (nonplanted) Galton--Watson tree $\t'$. Next,
graft the root of $\t'$ at the pointed vertex $v_h$ of $\t_*^{v_h}$
to form a
planted tree, which we denote by $\bar\t''$. Then $(\bar\t'',v_h)$
has the ``distribution'' $\mathrm{GW}_*$.

Otherwise said, the ``law'' under
$\mathrm{GW}_*$ of the distance $|v|$ of the pointed vertex $v$ to the
base is the counting measure on $\N$, and conditionally on $|v|=h$,
the subtrees $\bar T^{v}$ and $T_v$ are independent. More precisely, $T_v$
is a usual,\vspace*{1pt} nonplanted Galton--Watson tree with offspring
distribution $\nu$, and $\bar T^{v}$ has the law $\mathrm{GW}_*^h$.

It is then easy to see that for every $h\geq1$, re-planting the tree
$\t^{v_h}_*$ at the
vertex $v_h$ leaves its distribution $\mathrm{GW}_*^h$ invariant, the
pointed vertex $\hat v_h$ in
the re-planted tree being the base vertex of
$\t^{v_h}_*$
(this tree is $\hat{\t}^{\hat{v}_h}_*$
in our notation). Indeed, the offspring of the vertices along the branch
$(v_0,v_1,\ldots,v_h)$ are the same in both trees, while
the subtrees pending from the different offspring of
$v_1,\ldots,v_{h-1}$ are left unchanged. We deduce that
$(\t_*^{v_h},\t')$ and $(\hat{\t}_*^{\hat v_h},\t')$ have same
``distribution,'' recalling that $h$ is
not a random variable, but
rather chosen according to the counting measure on $\N$. This
entails Proposition~\ref{secverif-hypoth-eqrefhy}.~%
\end{pf}

We now turn our attention to $n\mu_{n,\xi_n}(t)$,
the number of edges in the component containing the randomly picked
edge $\xi_n$ in the forest $\bar\t_n(t)$. Recall that the latter
results from deleting every edge with probability
$1-\exp(-t/\sqrt{n})$, independently of the other edges,
in the planted Galton--Watson tree $\bar\t_n$ conditioned to have $n$ edges.
In this direction, it is convenient to introduce the following notation.
If $T$ is a rooted tree and $k\geq0$, we write $Z_k(T)$ for the
number of vertices at generation $k\geq
0$ in $T$, that is, at distance $k$ from the root. If the
tree is planted, then the definition of $Z_k(\bar T)$ is similar, but
counting only the vertices distinct from the base.
%
%
\begin{corollary}\label{C2} In the preceding notation, we have
\[
\E\bigl(n\mu_{n,\xi_n}(t)\bigr) \leq
\exp(-t/\sqrt{n})+2\sum_{k\geq
1}^{\infty}\exp(-kt/\sqrt{n}) \sup_{m\geq
1}\mathbb{E}\bigl(Z_k(\mathcal{T}_m)\bigr)
\]
\end{corollary}
\begin{pf}
%
For a vertex $u\in V(\mathcal{T}_n)$, let $e_u$ be the edge pointing
down from $u$ to the base, and for an edge $e$ of $\mathcal{T}_n$, let
$v(e)$ be the vertex such that $e_{v(e)}=e$. Let also $d(u,v)$ be the
graph distance in $\mathcal{T}_n$ between the vertices $u,v\in
V(\mathcal{T}_n)$.

Observe first that for every vertex $u\in V(\mathcal{T}_n)$, the edge
$e_u$ counts in the enumeration of $n\mu_{n,\xi_n}(t)$ if and only if
no edge on the path from $e_u$ to $\xi_n$ has been removed at time $t$.
Conditionally given $\mathcal{T}_n, \xi_n$, and for a given vertex
$u\in V(\mathcal{T}_n)$, this happens with probability
$\exp(-(d(u,v(\xi_n))+1)t/\sqrt{n})$ if $v(\xi_n)$ is an ancestor of
$u$, and with probability $\exp(-d(u,v(\xi_n))t/\sqrt{n})$ otherwise.
By distinguishing the vertex $v(\xi_n)$, for which the first formula
holds, from the other vertices, we thus have
\begin{eqnarray*}
\mathbb{E}\bigl(n\mu_{n,\xi_n}(t)\bigr)&\leq&
\mathrm{e}^{-t/\sqrt{n}}+\mathbb{E} \biggl(\sum_{u\in
V(\mathcal{T}_n)\setminus\{v(\xi_n)\}}
\mathrm{e}^{-d(u,v(\xi_n))t/\sqrt{n}}\biggr)\\
&=& \mathrm{e}^{-t/\sqrt{n}}+\frac{1}{n}\mathbb{E} \biggl(\sum_{u,v\in
V(\mathcal{T}_n),u\neq
v}\mathrm{e}^{-d(u,v)t/\sqrt{n}}\biggr),
\end{eqnarray*}
where the second identity follows from the fact that given
$\mathcal{T}_n$, $v(\xi_n)$ has the uniform law in $V(\mathcal{T}_n)$.

Next, notice that the set of pointed trees $(\bar T,v)$ with exactly
$n$ edges has $\mathrm{GW}_*$-measure equal to $n\P(|\t|=n)$, a
quantity which is strictly positive and finite by hypothesis.
So the conditional law $\mathrm{GW}_*(\cdot\mid|V(T)|=n)$ on the
space of
pointed trees with $n$ edges is well defined, and corresponds to the
distribution of $(\bar\t_n, \eta)$
where given $\t_n$, $\eta$ is a uniformly chosen vertex in $V(\t_n)$.

Combining these observations, we deduce that
\[
\E\bigl(n\mu_{n,\xi_n}(t)\bigr)=
\mathrm{e}^{-t/\sqrt{n}}+\mathrm{GW}_*\biggl(\sum_{u\in V(T)\setminus\{v\}}
\mathrm{e}^{-d(u,v)t/\sqrt{n}}\Bigm|\bigl|V(T)\bigr|=n\biggr).
\]
By definition, the number of vertices $u\in V(T)$ at distance $k\geq1$
from the pointed vertex $v$ equals
$Z_{k-1}(\hat{T}^{\hat{v}})+Z_k(T_v)$. 
Therefore,
%
%
\begin{eqnarray*}
\sum_{u\in V(T)\setminus\{v\}}\mathrm{e}^{-d(u,v)t/\sqrt{n}}&=&\sum_{k\geq
1}\mathrm{e}^{-kt/\sqrt{n}}\bigl(Z_{k-1}\bigl(\hat{T}^{\hat{v}}\bigr)+Z_k(T_v)\bigr)\\
&\leq&\sum_{k\geq 1}\mathrm{e}^{-kt/\sqrt{n}}\bigl(Z_k\bigl(\hat{T}^{\hat{v}}\bigr)
+Z_k(T_v)\bigr),
\end{eqnarray*}
where in the second step we performed a change of index. We conclude
that
\[
\mathbb{E}\bigl(n\mu_{n,\xi_n}(t)\bigr)\leq \mathrm{e}^{-t/\sqrt{n}}+\sum_{k\geq
1}\mathrm{e}^{-kt/\sqrt{n}} \mathrm{GW}_*
\bigl(Z_k\bigl(\hat{T}^{\hat{v}}\bigr)+Z_k(T_v) \mid  \bigl|V(T)\bigr|=n\bigr).
\]
By Proposition~\ref{secverif-hypoth-eqrefhy}, we have on the one hand
that
\begin{eqnarray*}
\mathrm{GW}_* \bigl(Z_k\bigl(\hat{T}^{\hat{v}}\bigr)
\bigm|\bigl|V(T)\bigr|=n \bigr) &=&\mathrm{GW}_* \bigl(Z_k\bigl( T^v
\bigr)\bigm|\bigl|V(T)\bigr|=n \bigr)
\\
&\leq&\mathrm{GW}_*\bigl(Z_k(T) \mid \bigl|V(T)\bigr|=n\bigr)
\\
&=&\mathrm{GW}\bigl(Z_k(T)\mid\bigl|V(T)\bigr|=n\bigr)
\\
&=&\E\bigl(Z_k(\t_n)\bigr).
\end{eqnarray*}
On the other hand, we saw in the proof of Proposition
\ref{secverif-hypoth-eqrefhy} that under $\mathrm{GW}_*$, the trees
$T^v$ and $T_v$ are independent, with $T_v$ having law
$\mathrm{GW}$. Therefore, since $|V(T)|=|V(T^v)|+|V(T_v)|-1$, we have
\begin{eqnarray*}
&& \mathrm{GW}_* \bigl(Z_k(T_v)\mid\bigl|V(T)\bigr|=n \bigr)
\\
&&\qquad=\sum_{m=1}^{n}\mathrm{GW}_*
\bigl(Z_k( T_v)\mid\bigl|V(T_v)\bigr|=m,\bigl|V
\bigl(T^v\bigr)\bigr|=n-m+1 \bigr)\\
&&\hspace*{16pt}\qquad\quad{}\times\mathrm{GW}_*\bigl(\bigl|V(T_v)\bigr|=m
\mid\bigl|V(T)\bigr|=n\bigr)
\\
&&\qquad=\sum_{m=1}^{n}\mathrm{GW}
\bigl(Z_k( T)\mid\bigl|V(T)\bigr|=m \bigr)\mathrm{GW}_*\bigl(\bigl|V(T_v)\bigr|=m
\mid\bigl|V(T)\bigr|=n\bigr)
\\
&&\qquad\leq \sup_{m\geq1}\E\bigl(Z_k(\t_m)\bigr)\sum
_{m\geq
1}\mathrm{GW}_*\bigl(\bigl|V(T_v)\bigr|=m
\mid\bigl|V(T)\bigr|=n\bigr)
\\
&&\qquad=\sup_{m\geq1}\E\bigl(Z_k(\t_m)\bigr),
\end{eqnarray*}
which completes the proof.
\end{pf}

Lemma~\ref{lem3} now follows readily from the following result by
Janson (see Theorem~1.13 in~\cite{Janson});
there exists some finite constant $C''$ depending only on the offspring
distribution $\nu$, such that
\[
\sup_{m\geq1}\E\bigl(Z_k(\t_m)\bigr)\leq C''k
\]
for every $k\geq1$. Indeed, we derive from Corollary~\ref{C2}
that for every $n\geq2$,
\[
\mathbb{E}\bigl(\mu_{n,\xi_n}(t)\bigr)\leq
\frac{\mathrm{e}^{-t/\sqrt{n}}}{n}+\frac{2C''}{n} \sum_{k\geq
1}k\mathrm{e}^{-kt/\sqrt{n}} \leq  \frac{ C\exp(-t/\sqrt
n)}{n(1-\exp(-t/\sqrt n))^2}.
\]

\section*{Acknowledgments}
G. Miermont acknowledges support and hospitality of CNRS/PIMS UMI 3069 at UBC,
Vancouver, where part of this research was done. We thank Daphn{\'e}
Dieuleveut for her careful reading of this paper.


%

\printaddresses

\end{document}